\documentclass[]{interact}

\usepackage{amsfonts,mathtools,xcolor}
\usepackage[round]{natbib}

\def\NAT@def@citea{\def\@citea{\NAT@separator}}
\makeatother

\theoremstyle{plain}
\newtheorem{theorem}{Theorem}[section]

\newtheorem{proposition}[theorem]{Proposition}

\theoremstyle{definition}

\theoremstyle{remark}

\newcommand{\re}{\mathbb{R}}
\newcommand{\N}{\mathbb{N}}
\newcommand{\ponto}{\, \cdot  \, }
\newcommand{\formu}[1]{\begin{eqnarray}\label{#1}}
\newcommand{\formub}{\begin{eqnarray} \nonumber}
\newcommand{\eformu}{\end{eqnarray}}

\newcommand{\aaa}{{\cal A}}

\newcommand{\xx}{\mathbf {x}}

\newcommand{\ppp}{{\cal P}}
\newcommand{\xxx}{{\cal X}}
\newcommand{\ie}{{\it i.e. }}
\newcommand{\Proof}{\vskip0.4cm \noindent {\bf Proof: }}
\newcommand{\eproof}{ \mbox{}\hfill$\sqcup\!\!\!\!\sqcap$ \vskip0.4cm \noindent}
\newcommand{\deffeq}{\mathrel{\overset{\makebox[0pt]{\mbox{\normalfont\tiny\sffamily def}}}{=}}}

\begin{document}
                               
\title{On the Bias of the Score Function of Finite Mixture Models}
\author{
\name{R. Labouriau\textsuperscript{a}\thanks{Correspondence to R. Labouriau. Email: rodrigo.labouriau@math.au.dk}}
\thanks{Published in: Labouriau, R. (2023). \emph{Communications in Statistics: Theory and Methods.} vol.{\bf 52}, issue {\bf 13}, pp 4461-4467.}
\affil{\textsuperscript{a} Department of Mathematics, Aarhus University. Ny Munkegade 118, Aarhus, 8000 C. DK}
}

\maketitle

\begin{abstract}
\noindent
We characterise the unbiasedness of the score function, viewed as an inference function for a class of finite mixture models. The models studied represent the situation where there is a stratification of the observations in a finite number of groups. 
We show that, under mild regularity conditions, the score function for estimating the parameters identifying each group's distribution is unbiased. We also show that if one introduces a mixture in the scenario described above so that for some observations, it is only known that they belong to some of the groups with a probability not in $\{ 0, 1 \}$, then the score function becomes biased. We argue then that under further mild regularity, the maximum likelihood estimate is not consistent. 
The results above are extended to regular models containing arbitrary nuisance parameters, including semiparametric models.
\end{abstract}

\section{Introduction}

It is well known that the score function is an unbiased estimating function under mild regularity conditions on the underlying statistical model  \citep[see][]{Godambe1960, GodambeThompson1976}. This fundamental property is used in the classic proofs of consistency of the maximum likelihood estimate \citep[see][]{JorgensenLabouriau2012,BarndorffNielsen1988}. 
We will present a collection of examples based on finite mixture models for which the score function is a biased inference function, and the maximum likelihood estimate for those models is not consistent. The proof presented is simple and involves relatively elemental mathematical tools.

The type of finite mixture models discussed here naturally appears in many statistical contexts and applications. Informally, these models represent the situation where there is a stratification of the observations in a finite number of groups, say $K$ disjoint groups. 
Suppose that the observations belonging to the same group follow the same distribution and that the $K$ distributions associated with each group are distinct elements of a sufficiently regular parametric family of probability measures. We show that in this scenario, under mild regularity conditions, the score function for estimating the parameters identifying each group's distribution is unbiased. We also show that if one introduces a mixture in the scenario described above, the score function becomes biased. More precisely, suppose now that there are some observations for which we do not know which group they belong to, still, instead, we might attribute probabilities of those belonging to different groups (and not all probabilities are in $\{ 0, 1\}$). In that case, the score function for estimating the parameters identifying each group's distribution is \emph{biased}.

In section \ref{Section.01} we briefly review some basic definitions and properties of inference functions and in section \ref{Section.02} we define the mixture model described above formally. The result characterising when the score function is biased is proved in section \ref{Section.02} 
and extended to sufficiently regular models with nuisance parameters in Section \ref{Section.04}.

\section{Basic Facts on Inference Functions} \label{Section.01}

Consider a parametric family of probability measures 
\begin{eqnarray} \label{Eqn01}
 \ppp = \left \{  
                P_\theta : \, \theta \in \Theta \subseteq \re^p , \Theta \mbox{ open}
            \right \}
\end{eqnarray}
defined on a common measure space  $\left ( \xxx , \aaa, \nu \right )$ which is denoted the \emph{basic model}. We assume that $\ppp$ is dominated by a 
$\sigma$-additive measure $\nu$ and that suitable versions of the 
Radon-Nykodim  derivatives, 
$d P_\theta / d\nu ( \, \cdot \, ) = p ( \, \cdot \, ;\theta) $,
are chosen to represent each probability measure $P_\theta$ in $\ppp$.
Moreover, suppose that the basic model given by \eqref{Eqn01}
is sufficiently regular to allow standard (well behaved) likelihood-based inference on $\theta$.

A function 
$\Psi : \xxx \times \Theta \longrightarrow \re^p$
such that for each 
$\theta\in\Theta$ the component function $\Psi (\ponto ; \theta)$ is measurable is called an \emph{inference function}. 
Given an inference function $\Psi$ and a sample $\mathbf{x} = (x_1, \dots, x_n) $
we define an estimate $\hat\theta = \hat\theta (\mathbf{x})$ of the parameter $\theta$ as the solution to the equation
\begin{equation}\label{eq.1.00}
 \mathbf{\Psi_n}( \mathbf{x}; \hat\theta )  =  1/n \sum_{i= 1}^n \Psi (x_i; \hat\theta ) = 0 \, . 
 \end{equation}
 The inference function $\Psi$ is said to be an \emph{unbiased inference function} when
$$
\int_\xxx \Psi (x, \theta) P_\theta (dx) = 0, \,\, \forall \theta \in \Theta \, .
$$
The inference function $\Psi$ is a \emph{regular inference function} if it is unbiased and the properties (i) - (v) below hold for all $\theta \in \Theta$ and all $i,j$ in $\{ 1, \dots , p\}$:
\begin{itemize}
\item[(i)]  The partial derivative $\partial \Psi (x;\theta )/\partial
\theta _i$ exists for $\nu $-almost every $x\in {\cal X}$ ;
\item[(ii)]  The order of integration and differentiation may be
interchanged as follows: 
$$
\frac \partial {\partial \theta _i}\int_{{\cal X}}\Psi (x;\theta )p(x;\theta
)d\mu (x)=\int_{{\cal X}}\frac \partial {\partial \theta _i}\left[ \Psi
(x;\theta )p(x;\theta )\right] d \nu (x)\,\,; 
$$
\item[(iii)]  $\mbox{E}_\theta \{\psi _i(\theta )\psi _j(\theta )\}\in \re$ and the 
$p\times p$
matrix 
$
V_\psi (\theta )=\mbox{E}_\theta \{\Psi (\theta )\Psi ^{\top }(\theta )\} 
$
is positive-definite;
\item[(iv)]  $\mbox{E}_\theta \left\{ \frac{\partial \psi _i}{\partial \theta
_r}(\theta )\frac{\partial \psi _j}{\partial \theta _s}(\theta )\right\} \in \re$ and the 
$p\times p$
matrix 
$
S_\psi (\theta )=\mbox{E}_\theta \{\nabla _\theta \Psi (\theta )\}\,\, 
$
is nonsingular.
\end{itemize}
Here $\nabla _\theta$ is the $p$-dimensional differential operator.
It can be shown (see \cite{JorgensenLabouriau2012}, Chapter 4) that the estimate sequence $\{ \hat \theta _n \}$ is consistent and asymptotically normally distributed.

\section{Basic setup} \label{Section.02}

We consider below the situation where only $K$ (for a fixed $K\in \N$) distributions of $\ppp$ are realisable. More precisely, suppose that we observe the realisations of  $I$ ($I\in \N $)  independent random variables, say
 $X_1, \dots , X_I$.
Suppose that there exists a finite set 
$\Theta_K = \left \{\theta_1, \dots , \theta_K \right \} \subset \Theta$ such that for each $i\in \left \{ 1, \dots , I \right \}$, 
there is one, and only one, $\theta_i \in \Theta_K$ such that  $X_i \sim P_{\theta_i}$.
Here the parameters $\theta_1, \dots , \theta_K$ determine the probability laws of $K$ different groups of observations termed the $K$ {\it components} of the basic model $\ppp$. It is assumed that the distributions of the components, $P_{\theta_1}, \dots , P_{\theta_K}$, are all different. 

In the setup defined above, each observation belongs to one, and only one, component. 
We will consider below scenarios ranging from the situation where it is known to which component each observation belongs to (complete information or absence of mixture), to the situation where this information is only known for some or none of the observations (partial or complete mixture, respectively). In order to characterize these scenarios, we introduce a $I\times K$ matrix 
$ \Pi = \left [ \pi_{i \, k} \right ]_{i=1,\dots,I; k=1, \dots, K}$, such that
for $i=1,\dots,I$ and $k=1, \dots, K$,
\begin{eqnarray} \label{Eqn02}
 \pi_{i \, k} \in [0,1] 
 \mbox{ and } 
 \sum_{j=1}^K \pi_{i \, j} = 1 \, .
\end{eqnarray}
The {\it finite mixed model} defined by the basic model $\ppp$ and the mixing matrix $\Pi$ is the model that for $i =1, \dots , I$
atributes the probability law given by the density
\begin{eqnarray} \label{Eqn03}
  \sum_{j=1}^K \pi_{i \, j} p\left (\, \cdot \, ; \theta_k \right ) \, ,
\end{eqnarray}
to the  $i^{th}$ observation.
Here $\pi_{i \, k}$ plays the role of the probability that the $i^{th}$ observation belongs to the $k^{th}$ component. 
If $\pi_{i \, k} \in \{ 0 , 1 \}$ for all  $i\in \{1,\dots,I\}$ and all $k \in \{1, \dots, K\}$, then we are in the situation where it is known
to which component each observation belongs to, which is equivalent to have 
a statistical model with a classification factor 
(the explanatory variable indicating the component, e.g. 
if the distributions of the  components are normal distributions 
with different means and a common variance, 
then we have a classic one-way analysis of variance model).
The cases where $\pi_{i \, k} \in [ 0 , 1 )$ for all  (or some)  
$i\in \{1,\dots,I\}$ correspond to the situation where the 
components are unknown for all (or for some) observations 
and therefore there is complete (or partial) mixture; 
we say then that the model given by (\ref{Eqn03}) {\it contains a mixture}.

\newpage
\section{Bias of the score function of finite mixed models}\label{Section.03}

The proposition below shows that score function is biased when a finite mixed model contains a mixture. This will imply (under further regularity conditions) that the maximum likelihood estimate for finite mixed model containing a mixture is not consistent.
\begin{proposition} \label{Theory.1}
Suppose that the statistical model related to the family $\ppp$ is sufficiently 
regular so that the score function 
$
\mathbf {S}(\xx, \theta)=\sum_{i=1}^I S_i( x_i, \theta)
$ 
is well defined and unbiased, i.e., for all $\theta \in \Theta$ and all $i\in \{1, \dots , I\}$, 
 \begin{eqnarray}\label{Eqn04}
 E_\theta \left [ S_i(X ,\theta) \right ] = 
 \int_\xxx S_i( x; \theta) p(x; \theta)  d \nu (x) = 0 
 \, .
\end{eqnarray}
Then the score function of the mixed model defined by the basic model 
$\ppp$ and the mixing matrix $\Pi$ satisfying (\ref{Eqn02}) is unbiased 
if, and only if, 
the finite mixed model 
does not contain a mixture. 
\end{proposition}
\Proof 

($\Longrightarrow$) 
If the model does not contain a mixture, then  $\pi_{i \, k} \in \{ 0 , 1 \}$ 
for all  $i\in \{1,\dots,I\}$ and all $k \in \{1, \dots, K\}$.
Therefore the expectation of the component of the score function 
with respect to the parameter $\theta_k$ ($k=1, \dots ,K$) is given by
\begin{eqnarray}\nonumber
  E_{\theta_k} \left [ S_i(X_i; \theta_k) \right ] = 
 \int_\xxx \pi_{i \, k} S_i(x; \theta_k) p(x; \theta_k)  d \nu (x) =
 \mbox{( by (\ref{Eqn04}) )}  = 0 
\end{eqnarray}
and therefore the score function is unbiased.

($\Longleftarrow$)
Suppose that the model given by (\ref{Eqn03}) contains a mixture. 
We proof that in this case the score function is biased. 
Since the model given by (\ref{Eqn03}) contains a mixture, then there exist an $i \in \{1, \dots ,I\}$ and 
a $k$ and a $k^\prime \in \{1, \dots ,K \}$ ($k \ne k^\prime$)
such that $0 < \pi_{i \, k} < 1$ and $0 < \pi_{i \, k^\prime} < 1$. 
Assume, without loss of generality, that $\pi_{i \, j}=0$ for all 
$j \in \{1,  \dots ,K \} \setminus \{k, k^\prime \}$.
Denoting the differentiation operator with respect to $\theta_k$ by $\Delta_{\theta_k}$, we have that
\begin{eqnarray}\label{Eqn05}
  E_{\theta_k} \left [ S_i(X_i; \theta_k) \right ] 
  & = &
          \int_\xxx 
          \frac{
              \Delta_{\theta_k} 
              \sum_{j=1}^{K} \pi_{i \, j} p(x; \theta_j)
           }
           {
           \sum_{j=1}^{K} \pi_{i \, j}p(x; \theta_j)
           }
   p(x; \theta_k)  d \nu (x)
   \\ \nonumber & = &
   \int_\xxx 
          \frac{
              \pi_{i \, k}
              \Delta_{\theta_k} 
               p(x; \theta_k)
           }
           {
            \pi_{i \, k}p (x; \theta_k) + \pi_{i \, k^\prime}p (x; \theta_k^\prime)
           }
           p(x; \theta_k)  d \nu (x)
  \\ \nonumber & < &
   \int_\xxx 
          \frac{
              \pi_{i \, k}
              \Delta_{\theta_k} 
               p(x; \theta_k)
           }
           {
            \pi_{i \, k}p (x; \theta_k) 
           }
           p(x; \theta_k)  d \nu (x)
     \\ \nonumber & &  \mbox{(since $0<\pi_{i \, k}<1$)}
     \\ \nonumber & < &
 \int_\xxx S_i(x; \theta_k) p(x; \theta_k)  d \nu (x) =
 \mbox{( by (\ref{Eqn04}) )}  = 0 \, ,
\end{eqnarray}
therefore the score function is biased.
\eproof

 
\vspace{0.2cm}
\noindent
\emph{Remark 1:} 
It is well known that the score function is unbiased under mild regularity 
conditions on the statistical model \citep[see][]{JorgensenLabouriau2012}; 
therefore, it is natural to ask where the proof of unbiasedness of the score 
function fails in the present example. 
The following sketch of the proof of unbiasedness of the score function
enlighten this point. Denoting the density of an element of the basic model $\ppp$ by $f(\,\cdot\, ; \theta)$, we have, for each $\theta\in\Theta$,
\begin{eqnarray} \nonumber 
  E_{\theta} \left [ S(X; \theta) \right ] 
  & = & 
          \int
          S(x; \theta)  f(x; \theta)  d \nu (x)
      = \int \frac{\partial}{\partial \theta} \log \{ f(x;\theta) \}  f(x; \theta)  d \nu (x)
   \\ \label{Eqn06} & = &  \!\!\!\!
          \int \frac{\partial}{\partial \theta} f(x;\theta) 
                  \frac{1}{f(x; \theta)}  f(x; \theta)  d \nu (x)
          \,  =
           \int \frac{\partial}{\partial \theta} f(x;\theta) d \nu (x) 
   \\ \label{Eqn07} & = &    \!\!    
       \frac{\partial}{\partial \theta} \int f(x;\theta) d \nu (x)  
     \,  =  \frac{\partial}{\partial \theta} 1 \, = 0 \, .
\end{eqnarray}
Here the crucial part of the argument is the differentiation under the integration sign (a property obtained for example by assuming that the tails of the density function $f(\, \cdot \,; \theta )$ decays sufficiently rapidly) used in (\ref{Eqn07}), but this condition typically holds for all densities of the family $\ppp$ and therefore it holds also for mixtures of densities of $\ppp$; therefore the classic proof typically does not fail at this point. Indeed, the classic proof of unbiasedness of $S$ breaks before this point, namely at the factorization step used in (\ref{Eqn06}). 

\vspace{0.3cm}
\noindent
\emph{Remark 2:} 
\noindent
We sketch below an argument showing that, under mild regularity conditions, the maximum likelihood estimate is not consistent when the score function is biased. 
Since the purpose here is just to be illustrative, we consider the case where $\Theta \subseteq \re$ and we assume all the regularity conditions required to keep the argument simple. Here we use the notion of consistency in probability \citep[see][page 232]{Lehmann1983}, \ie, a sequence of estimates $\{\hat\theta_n\}$ is consistent if, and only if,
$\hat\theta_n \xrightarrow[n\rightarrow \infty]{P_\theta} \theta $, for all $\theta\in\Theta$. 

Suppose, by hypothesis of absurdum, that $S$ is biased and that there is a sequence $\{\hat\theta_n\}$ of roots of $1/n \sum_{i= 1}^n S (X_i; \theta )$ that is consistent.
Define for each $\theta_*\in\Theta$ the function $\lambda_{\theta_*} : \Theta \longrightarrow \re$ given by 
$\lambda_{\theta_*} (\theta ) = \int_\xxx S(x; \theta ) p(x; \theta_*)  d \nu (x)$ and assume that the function $\lambda_{\theta_*} $ is continuous for any choice of $\theta_*\in\Theta$. 
Clearly, the score function $S$ is unbiased if, and only if, 
$\lambda_{\theta} (\theta) = 0 $ for all $\theta\in\Theta$. 
If $S$ is biased, then there exist $\theta_*\in\Theta$ such that $\lambda_{\theta_*} (\theta_*) \ne 0 $ and, since  $\lambda_{\theta_*}$ is continuous in a neighbourhood of $\theta_*$,  there exist an $\epsilon = \epsilon_{\theta_*}  >0$ such that for all $\theta \in (\theta_* - \epsilon, \theta_* + \epsilon )$, $\lambda_{\theta_*} (\theta ) \ne 0$. Since  $
 1/n \sum_{i= 1}^n S (X_i; \theta_* ) \xrightarrow{P_{\theta_*-a.s.}} \lambda_{\theta_*} (\theta_* )$, then for $n$ sufficiently large  $1/n \sum_{i= 1}^n S (X_i; \theta )$ will not have a root in $(\theta_* - \epsilon, \theta_* + \epsilon )$. Therefore, for $n$ sufficiently large, $\theta_n - \theta_* > \epsilon$, $P_{\theta_*}-a.s.$, which contradicts the assumption of consistency of  $\{\hat\theta_n\}$.

\section{
Extension to models with nuisance parameters
} \label{Section.04}
It is of interest to extend the results above for models containing nuisance parameters (as suggested by a reviewer). We provide some details below. 
Consider the family of probability measures 
\begin{align*}
 \widetilde\ppp = \left \{  
                \widetilde P_{\theta \lambda}: \, \theta \in \Theta \subseteq \re^p , \Theta \mbox{ open},
                \lambda\in\Lambda
            \right \}
\end{align*}
defined on a common measure space  $\left ( \xxx , \aaa, \nu \right )$.
Here $ \widetilde\ppp$ will play the role of the basic model as in the construction of a finite mixed model given in Section \ref{Section.02}.
Moreover, $\theta$ is a parameter of interest and $\lambda$ a nuisance parameter of arbitrary nature. 
As before, we assume that $\widetilde\ppp$ is dominated by $\nu$ and that suitable versions of the Radon-Nykodim  derivatives, 
$d \widetilde P_{\theta \lambda}/ d\nu ( \, \cdot \, ) =  \tilde p ( \, \cdot \, ;\theta, \lambda) $,
are chosen to represent each probability measure $\widetilde  P_{\theta \lambda}$ in $\widetilde\ppp$.
We assume that the parametrisation of $\widetilde\ppp$ is identifiable, \ie, 
the map 
$(\theta, \lambda) \mapsto \widetilde P_{\theta\lambda}$ 
from $\Theta\times\Lambda$ to $\widetilde\ppp$ is a bijection, and that
$\widetilde\ppp$ is sufficiently regular to allow 
likelihood-based inference on $\theta$ in the presence of the nuisance parameter $\lambda$.

A function $\Psi : \xxx \times \Theta \times \Lambda \longrightarrow \re^q$ such that for each 
$(\theta, \lambda) \in\Theta\times\Lambda$ the component function $\Psi (\ponto ; \theta, \lambda)$ is measurable is called a \emph{quasi-inference function} 
(see \citeauthor{JorgensenLabouriau2012}, \citeyear {JorgensenLabouriau2012}, Chapter 4, section 4.4 for details on inference functions in the presence of nuisance parameters). 
In this context, an \emph{inference function} is a quasi-inference function that does not depend on the nuisance parameter. A quasi-inference function $\Psi: \xxx \times \Theta \times \Lambda \longrightarrow \re^q$ is regular when the conditions (i) - (v) given in Section  \ref{Section.01}, with $p(\,\cdot\, , \theta)$ replaced by $\tilde p(\,\cdot\, , \theta, \lambda)$,  hold for all $\theta \in \Theta$, all $\lambda\in\Lambda$, and all $i,j$ in $\{ 1, \dots , q\}$.
In particular, $\Psi$ is unbiased when
\begin{align*}
\int_\xxx \Psi (x, \theta, \lambda) \widetilde P_{\theta\lambda} (dx) = 0,
\mbox{ for all } \theta \in \Theta \mbox{ and all } \lambda\in\Lambda \, .
\end{align*}
We assume that the partial score function for estimating $\theta$ in $\widetilde\ppp$ given by 
\begin{align}\label{eq:partialscore}
 \mathbf {\widetilde S}(\xx \, , \theta , \lambda) = 
 \sum_{i=1}^I \nabla_\theta \log \tilde p( x_i , \theta , \lambda), \mbox{for } (\theta, \lambda) \in\Theta\times\Lambda ) , 
 \end{align}
 is well defined.

In the context of models with nuisance parameter, the \emph{finite mixed model defined by the basic model $\widetilde\ppp$ and the mixing matrix $\Pi$} is the model that attributes to the $i^{th}$ observation (for $i =1, \dots , I$) the probability law given by the density
\begin{eqnarray} \label{Eqn03b}
 \sum_{j=1}^K \pi_{i \, j} \tilde p\left (\, \cdot \, ; \theta_k , \lambda_0 \right ) \, ,
\end{eqnarray}
for some  $\lambda_0\in\Lambda$.
Here $\lambda_0$ is an unknown but fixed point of $\Lambda$.
The following proposition, analogue to Proposition \ref{Theory.1}, holds in the context consider here.
\begin{proposition} \label{Theory.2}
Suppose that the statistical model related to the family $\widetilde\ppp$ is sufficiently 
regular so that the partial score function  for estimating $\theta$ under $\widetilde\ppp$,
given in \eqref{eq:partialscore}, is well defined and unbiased.
Then the score function for estimating $\theta$ under the finite mixed model defined by 
$\widetilde\ppp$ and the mixing matrix $\Pi$, satisfying (\ref{Eqn02}),  is unbiased 
if, and only if, the finite mixed model defined by $\widetilde\ppp$
does not contain a mixture. 
\end{proposition}

\noindent
The following definition is required for the proof. For each $\lambda$ in $\Lambda$, define
the \emph{$\lambda$-orbit of $\widetilde\ppp$} by
\begin{align*}
\widetilde\ppp_{\lambda} \deffeq
\left\{ 
 P_{\theta \lambda} \in \widetilde\ppp : \theta\in\Theta \, .
\right\}
\end{align*}
Clearly, 
$
\widetilde\ppp = \cup_{\lambda\in\Lambda}  \widetilde\ppp_{\lambda} 
$.
Moreover, a quasi-inference function $\Psi$ is unbiased if, and only if, for each $\lambda\in\Lambda$ the restriction of $\Psi$ to $\widetilde\ppp_{\lambda}$ is an unbiased inference function. 
\Proof 
The finite mixed model defined by 
$\widetilde\ppp$ and the mixing matrix $\Pi$ contains a mixture if, and only if, 
the finite mixed model defined by the $\lambda_0$-orbit $\widetilde\ppp_{\lambda_0}$ and the mixing matrix $\Pi$ contains a mixture. The proof will follow by applying the proposition \ref{Theory.1} to the restrictions of the score function to the orbit $\widetilde\ppp_{\lambda_0}$, as given below.

($\Longrightarrow$) 
 If the finite mixed model does not contain a mixture, then 
 $\widetilde\ppp_{\lambda_0}$ do not contain a mixture. Therefore, by the  proposition \ref{Theory.1}, the restriction of the score function to $\widetilde\ppp_{\lambda_0}$ is unbiased, implying that the score function for estimating $\theta$ under $\widetilde\ppp$ is unbiased.

($\Longleftarrow$) 
If the finite mixed model contains a mixture, then 
$\widetilde\ppp_{\lambda_0}$ contains a mixture. Therefore,  by the  proposition \ref{Theory.1}, the restriction of the score function to $\widetilde\ppp_{\lambda_0}$ is biased, implying that the score function for estimating $\theta$ under $\widetilde\ppp$ is biased.
\eproof




\begin{thebibliography}{plainnat}

\bibitem[\protect\citeauthoryear
{Barndorff-Nielsen}
{Barndorff-Nielsen}{1998}]  
{BarndorffNielsen1988}
Barndorff-Nielsen, O.E. (1988). {\it Parametric Statistical
Models and Likelihood\/}. Lecture Notes in Statistics, Springer-Verlag, New
York.

\bibitem[\protect\citeauthoryear
{Godambe} 
{Godambe}{1960}]
{Godambe1960} 
Godambe, V.P. (1960). An optimum property of regular maximum
likelihood estimation. {\it Ann. Math. Statist.\/} {\bf 81}, 1208--1212.

\bibitem[\protect\citeauthoryear
{Godambe and Thompson} 
{Godambe and Thompson}{1976}]
{GodambeThompson1976}
Godambe, V.P. and Thompson, M.E. (1976). Some aspects of the theory
of estimating equations. {\it J. Statist. Plann. Inference\/} {\bf 2},
95--104.

\bibitem[\protect\citeauthoryear
{J{\o}rgensen and Labouriau}
{J{\o}rgensen and Labouriau}{2012}]
{JorgensenLabouriau2012}
J\o{}rgensen, B. and Labouriau, R. (2012).
{\it Exponential Families and Theoretical Inference.\/}
{\em Monografias de Matematica \/}~{\em 52},  196p.
Instituto de Matematica Pura e Aplicada (IMPA), Rio de Janeiro, Brazil.
(https://pure.au.dk/portal/files/51499534/Mon52.pdf )

\bibitem[\protect\citeauthoryear
{Lehmann}
{Lehmann}{1983}]
{Lehmann1983}
Lehmann, E.L. (1983).
{\it Theory of Point estimation.\/}
John Wiley and Sons, New York.  506p.

\end{thebibliography}
\end{document}